# Languages, Algorithms, Procedures, Calculi and Metalogic

**Mark Burgin**

Department of Mathematics
University of California, Los Angeles
405 Hilgard Ave.
Los Angeles, CA 90095

**Abstract**.

Logicians study and apply a multiplicity of various logical systems. Consequently, there is necessity to build foundations and common grounds for all these systems. This is done in metalogic. Like metamathematics studies formalized mathematical theories, metalogic studies theories in logic, or logics. The discipline of logic has been developed with the aim to model and study human thinking and reasoning. A more realistic understanding relates logic only to reasoning. Reasoning is a mental and verbal activity. Any activity is based on actions and operations organized by procedures and algorithms. That is why procedures and algorithms are basic for the development of different logics, their study and application. In this work, we study structures used in logic with the aim to reflect goals of logic as a discipline. The main structure is a logical calculus, which pivotal for the contemporary logic. Logical calculus is considered here as a two-tier construction with a logical language on the first level and the calculus itself on the second level. A system of interdependencies that exists between properties of deductive calculi and utilized by these calculi algorithms are explicated and studied. Here only syntactic parts of logics, namely, deductive logical calculi and corresponding languages are considered. Semantic and pragmatic parts are studied elsewhere.

**Key words.** logic, algorithm, procedure, language, syntactic calculus, logical language, inference

# 1. Introduction

Languages are created and used for communication. Procedures and algorithms are created and used for action and transformation. Calculi combine languages, procedures and algorithms to provide means for derivation and generation of new entities from existing ones.

The word *calculus* has two meanings in mathematics. The most popular in general mathematics understanding is that *Calculus* is a name that is now used to denote the field of mathematics that studies properties of functions, curves, and surfaces. As this is the most popular meaning in mathematics, we call it *the calculus*. It is usually subdivided into two parts: *differential calculus* and *integral calculus*. The main tool of the calculus is operating with functions to study properties of these functions. This operation can be regarded as a generalized calculation with these functions. This explains the name *calculus* used for this field, which originated from the Latin word meaning pebble because people many years ago used pebbles to count and do arithmetical calculations. The Romans used *calculos subducere* for "to calculate."

Thus, the calculus is called so because it provides analytic, algebra-like techniques, or means of computing, which apply algorithmically to various functions and curves. Many mathematical problems that had very hard solutions or even such problems that mathematicians had not been able to solve, after the calculus had been developed, became easily solvable by mathematics students.

Later the calculus developed into analysis, or mathematical analysis. There are also other calculi in analysis, for instance, operational calculus and calculus of variations.

Another mathematical meaning of the word *calculus* comes from mathematical logic where *calculus* is a formal system used for logical modeling of mathematical and scientific theories. A logical calculus consists of three parts: axioms, rules of deduction (inference), and theorems (cf. (Kleene, 2002; Mendelson, 1997).

The idea of the concept of logical calculus comes from Leibniz, who also introduced names *differential calculus* and *integral calculus*. He wrote that in future informal and vague arguments of philosophers would be changes for formal and exact calculations with formulas. Such calculations would allow one to find who of those philosophers was right and who was wrong.

To make such generalized calculations with formulas, people use definite rules. Systems of rules form algorithms when they are precise, exactly realizable and sufficiently simple to be performed by a mechanical device. Otherwise, such systems are called procedures.

Generalized calculations are performed with symbolic expressions, which are elements of definite languages, usually, formal languages. The goal of this work is to study relations between algorithms, procedures, languages, and calculi as a part of metalogic.

Here we are interested mostly in logical languages and calculi, as well as in algorithms of deduction and formal inference. Algorithms form a foundation for logic as logic as field has been developed with the aim to model and study human thinking and reasoning. A more realistic approach relates logic only to reasoning. Reasoning is an activity. Any activity is based on actions and operations organized in procedures and algorithms. That is why algorithms are necessary for the development of different logics, their study and application.

Some think that algorithms are necessary only for syntax. However, semantics and pragmatics of any language and logic are also determined by corresponding algorithms given, as a rule, in a form of rules.

In this work, we study only syntactic parts of logics, namely, deductive logical calculi. Semantic and pragmatic parts are studied elsewhere. A logical calculus in the metalogical setting is constructed as a two-tier hierarchical system. The first level gives a constructive representation of a logical and algorithmic languages used in logic. This level serves as the base for the second level, which consists of a syntactic, or deductive, logical calculus. To build a mathematical theory (metatheory) of logic, we utilize constructions from the theory of logical varieties (Burgin, 1995; 2004).

Properties of constructive representations of logical languages (Section 2) and syntactic (deductive) logical calculi (Section 3) are studied.

It is necessary to remark that the algorithmic level of logic studied here form a foundation (metalogic) not only for a huge diversity of logics studied in mathematics, logic and computer science, but also for universal logic introduced by Béziau (1994).

**Denotations and basic definitions:**

$N$ is the set of all natural numbers. $\omega$ is the sequence of all natural numbers.

$\varnothing$ is the empty set. The logical symbol $\forall$ means "for any". The logical symbol $\exists$ means "there exists".

If $X$ is a set, then $P(X)$ is the set of all subsets and $P_{fin}(X)$ is the set of all finite subsets of $X$.

If $A$ is a system of algorithms and $X$ is a set, then $A(X)$ denotes the set of all elements that can be obtained by application of algorithms from $A$ to the elements from $X$. For instance, if algorithms from $A$ cannot be applied to the elements from $X$, then $A(X) = \varnothing$.

A *binary relation* $T$ between sets $X$ and $Y$ is a subset of the direct product $X \times Y$. The set $X$ is called the *domain* of $T$ and denoted by $D(T)$ and $Y$ is called the *codomain* of $T$ and denoted by $CD(T)$. The *range* of the relation $T$ is $R(T) = \{ y \,;\, \exists\, x \in X\, ((x, y) \in T)\}$. The *definability domain* of the relation $T$ is $DD(T) = \{ x \,;\, \exists\, y \in Y\, ((x, y) \in T) \}$.

An *n-ary relation* $Q$ in a set $X$ is a subset of the direct power $X^n$.

The sequential composition $T \circ P$ of a binary relation $T$ between sets $X$ and $Y$ and a binary relation $P$ between sets $Y$ and $Z$ is a subset of the direct product $X \times Z$ defined as $T \circ P = \{ (a, c)$ where $a \in X$, $c \in Z$, and there is $b \in Y$ such that $(a, b) \in T$ and $(b, c) \in P \}$.

The closure $T^*$ of a binary relation $T$ in a set $X$ is the union of all sequential powers of the binary relation $T$, i.e., $T^* = \bigcup_{n=1}^{\infty} T^n$.

A *function* or *total function* from $X$ to $Y$ is a binary relation between sets $X$ and $Y$ that satisfies two following conditions: 1) there are no elements from $X$ which are corresponded to more than one element from $Y$; and 2) some element from $Y$ is corresponded to any element from $X$. Often total functions are also called everywhere defined functions.

A *partial function f* from $X$ to $Y$ is a binary relation in which there are no elements from $X$ which are corresponded to more than one element from $Y$.

A word in an alphabet $X$ is any finite string of elements from $X$. The symbol $\varepsilon$ denotes the empty word. A formal language $L$ in an alphabet $X$ is any subset of the set $X^*$ of all words in the alphabet $X$.

### 2. Languages as Calculi

A definition of any logic starts with a definition of its language. A definition of any language starts with a definition of its alphabet. In the most general case, we assume that any set $X$ may be an alphabet of a language. This assumption is made to allow us to use, for example, a vocabulary as an alphabet and thus, to build a conventional (natural or artificial like programming) language from words and not only from symbols/letters. A formal language is any set of words in some alphabet. However, to really know a language, we need to have its more concrete definition. There are three forms of language definitions/representations: demonstrative, descriptive and constructive definitions.

**Definition 2.1.** A *demonstrative definition* of a language $L$ represents this language as a list (collection) of words.

When the language $L$ is not very big, such presentation can be complete. For infinite languages and even very big languages, this form allows only representation of a part of the language. For instance, decimal representation of natural numbers can be given in form of a list: 1, 2, 3, … or 1, 2, 3, … , n, n + 1, …

**Definition 2.2.** A *descriptive definition* of a language L gives some description of this language.

A conventional description of a language in the theory of formal languages is representation of this language by a formula.

*Example* 2.1. Any regular language is defined by a formula that is called a regular expression (Hopcroft, *et al*, 2001).

**Definition 2.3.** A *constructive definition* of a language L gives an algorithm (rules or operations) to construct this language.

In some cases, rules of language construction form a procedure.

Here we are mostly interested in constructive definitions of languages. The reason for this is that in its most general form, a logical language is usually treated as a set of rules for constructing formulas for some logic. Logic works with these formulas, deducing some formulas from others, transforming formulas, and assigning truth values to formulas based on the rules of that logic.

It is necessary to note that a constructive definition of a language L is a kind of a formal calculus (compare Definitions 2.4 and 3.1).

According to the axiomatic theory of algorithms (Burgin, 2004a), construction of a language can be organized in three main types:

1. *Production/computation* when an algorithm (system of algorithms) *A* (in general, potentially) builds all words from *L* and only such words.

2. *Acceptance* (or *separation*) when an algorithm (system of algorithms) *A* (in general, potentially) accepts all words from *L* and only such words.

3. *Decision* when an algorithm (system of algorithms) *A* (in general, potentially) accepts all words from *L* and rejects all other words.

Thus, we have three types of constructive language definitions/representations: production, acceptance, and decision definitions.

*Example* **2.2.** Context-free languages are usually defined by derivation of their words utilizing rules from a context-free grammar *G* (Hopcroft, *et al*, 2001). The grammar *G* is an algorithm that is used for production of a context-free language as derivation is a kind of production.

*Example* **2.3.** Context-free languages are also defined by recursive inference utilizing rules from a context-free grammar *G* (Hopcroft, *et al*, 2001). In this case, the grammar *G* is an algorithm that is used for acceptance of a context-free language as recursive inference is a kind of production.

**Remark 2.1.** To build languages, it is possible to use not only recursive inference and conventional derivation, but also inductive inference and inductive derivation, which are more powerful than recursive inference and conventional derivation (Burgin, 2005).

Constructive definition gives a more detailed representation of a language than two other kinds of representation. It is formalized by means of the corresponding fundamental triad (Burgin, 2004b) described in the following definition.

**Definition 2.4.** A language in a constructive representation/definition is a triad (a named set) of the form $L = (X, R, L)$ where $X$ is the alphabet, $R$ is the set of constructive algorithms/rules and $L$ is the set of words of the language *L*.

Logical languages are a special kind of artificial languages developed intentionally within a culture. The typical feature of logical languages is that their structure Inner relations) and grammar (formation rules) are intended to express the logical information within linguistic expressions in clear and effective ways. Languages used in logic have, as a rule, constructive definitions in a form of production rules. Elements of logical languages are logical expressions or formulas. To emphasize that these formulas are constructed in a proper way, they are often called well-formed formulas.

*Example* **2.4.** Elements of the language of the classical propositional or sentential logic/calculus give a formal representation of propositions. Propositional variables are denoted by the capital letters of the Latin alphabet (A, B, C, etc.) or the small letters of the Greek alphabet (χ, φ, ψ, etc.). However, only one alphabet is usually used, but the two are not mixed. These letters are considered as atomic formulas and form a part of the alphabet of the language. Another part is formed by the symbols denoting the following connectives (or logical operators): negation denoted by ⌐, logical "and" denoted by ∧, logical "or" denoted by ∨, implication denoted by →, and equivalence denoted by ↔. Logicians use other symbols to denote the same logical operators: negation is also denoted by ~ , logical "and" is also denoted by & and ·, implication is also denoted by ⇒ and ⊃, and equivalence is also denoted by ≡ and ⇔. It is possible to use fewer operators (and thus, a smaller alphabet) by expressing some of these operators by mean of others, e.g., P → Q is equivalent to ⌐P ∨ Q. For example, Church (1956) uses only one logical operator ⊃. In addition, the left and right parentheses, ( and ) or the left and right brackets [ and ] are included in the alphabet.

Elements of the language $L_P$ of the classical propositional or sentential logic/calculus are called well-formed formulas (*wffs*). To build the set of well-formed formulas (*wffs*) the following rules are used:

1. Letters of the alphabet are *wffs* from $L_P$.

2. If φ is a *wff*, then ⌐φ is a *wff* from $L_P$.

3. If φ and ψ are *wffs*, then (φ ∧ ψ), (φ ∨ ψ), (φ → ψ), and (φ ↔ ψ) are *wffs* from $L_P$.

These rules form the set of algorithms $R$ that build the language $L_P$ of the classical propositional calculus.

***Example* 2.5.** Elements of the language $L_{CPC}$ of the classical predicate logic/calculus of the first order give a formal representation of binary properties. The predicate calculus language has a developed alphabet, making heavy use of symbolic notation. Lower-case letters *a*, *b*, *c*, ..., *x*, *y*, *z*, ... are used to denote individuals. Upper-case letters *M*, *N*, *P*, *Q*, *R*, ... are used to denote predicates.

The alphabet $L_{CPC}$ of the language $L_{CPC}$ consists of the following parts:

- A set **F** of function symbols (common examples include + and · );

- A set **P** of predicate symbols (common examples include = and <);
- A set **C** of logical connectives (usually, it is ⌐, ∧, ∨, →, and ↔);
- A set **S** of punctuation symbols (usually, examples include ( , ), : and ,);
- A set **Q** of quantifiers (usuallly, they are ∀ and ∃ );
- A set **V** of variables.

Every function symbol, relation symbol, and connective is associated with an arity. Namely, an *n*-ary function has the form $f: X^n \to X$ and an *n*-ary function has the form $P(x_1, x_2, \ldots, x_n)$. The set of *n*-ary function symbols is denoted $\mathbf{F}_n$, and the set of *n*-ary predicate symbols is denoted $\mathbf{P}_n$. As a rule, 0-ary predicates and/or 0-ary functions are called constants. Another way to deal with constants is to include their names in the alphabet of the language.

The language $L_{CPC}$ of the classical predicate calculus encompasses the language $L_P$ of the classical propositional calculus as propositions may be formed by juxtaposition of a predicate with an individual.

Elements of the language $L_{CPC}$ of the classical predicate logic/calculus are also called well-formed formulas (*wff*s). To build the set of well-formed formulas (*wff*s) the following rules are used:

1. Letters of the alphabet are *wff*s from $L_{CPC}$.
2. If φ is a *wff*, then ⌐φ is a *wff* from $L_{CPC}$.
3. If φ and ψ are *wff*s, then (φ∧ψ), (φ ∨ ψ), (φ → ψ), and (φ ↔ ψ) are *wff*s from $L_{CPC}$.
4. If φ is a *wff* containing a free instance of variable *x*, then ∃*x*φ and ∀*x*φ are *wff*s from $L_{CPC}$.

Here a variable is free if it is not related to a quantifier. Consequently, transformation rule (4) makes any instance of *x* bound (that is, not free) in the formulas ∃*x*φ and ∀*x*φ .

***Example* 2.6.** Elements of the language $L_{QA}$ of the logic of questions and answers are questions and propositions (Belnap and Steel, 1976). Such logics are called erotetic. A formalized language of an erotetic logic consists of two parts: assertoric and erotetic. The assertoric part of $L_{QA}$ is usually a first-order language. As far as the assertoric part of $L_{QA}$ is concerned, the concepts of term, well-formed formulas, variables, deduction, etc., are defined in a conventional way. Questions are the meaningful expressions of the erotetic part of $L_{QA}$.

Usually questions are separated into relevant classes. One way of classifying questions is in terms of the surface characteristics that give such classes as:

1) yes/no questions (for example: "Is it now ten a.m.?");

2) item-specification questions (for example: "Who is a student?");

3) instruction-seeking questions (for example: "How to learn logic?").

Another way to classify questions is to take into account the nature of the answers. It gives us factual questions (for example: "What time is it now?"), normative questions (for example: "How it is necessary to drive when it is raining?"), and counter-factual questions (for example: "What might happen if we met a year ago?"). There are also special rules to build correct questions of the erotetic calculus.

Thus, we can see that logicians use a diversity of logical languages and continue to invent new ones.

Utilization of a logical language involves different operations with their elements. Such operations are performed according to definite rules (algorithms). The main operations are inference and substitution. Consequently, to build a logic, we need algorithms and it is natural to consider algorithmic languages, elements of which are texts/expressions that describe algorithms. Examples of algorithmic languages are the language of Turing machines, the language of finite automata or the language of inductive Turing machines (Burgin, 2005). All these algorithms can be used for inference. Although there is a difference between algorithms and their descriptions, as it is demonstrated in (Burgin, 2005), here we do not emphasize this difference and write, for simplicity, that algorithms belong to an algorithmic language.

Algorithmic languages have the same types of representations as logical languages or any languages have, that is, there are three classes of language representations: representative, descriptive, and constructive representations. As in a general case, a constructive representation of an algorithmic language has the form of a fundamental triad.

Systems of algorithms are often algebras or calculi because there are rules of composition of algorithms and of inference or derivation of algorithms in such a system. Examples of composition are sequential composition of two algorithms when the result of the first algorithm is given as the input for the second algorithms and parallel composition of two algorithms when the result of them composition algorithm consists of the results of both composed algorithms. An example of inference or derivation is given by the following reasoning: if in a monotone logical calculus **C,** there are inference algorithms *r* and *q*, then their sequential composition also belongs to this calculus.

However, there is an essential difference between logical and algorithmic languages with respect to problems of logic. In the syntactic context, semantics of a logical language, which is a language of formulas, can be ignored. For instance, it is possible to treat propositions simply as letters or words. In contrast to this, semantics and, especially, the dynamic semantics of algorithmic languages is pivotal for logic. The reason is that to use algorithms, we need rules telling us how to apply these algorithms. Such rules, or sometimes they are called metarules, form the dynamic semantics of algorithmic languages.

### 3. Logical Calculi

According to the theory of logical varieties (Burgin, 1997), there are three kinds of logical calculi or simply, logics:

1. *Deductive or syntactic* logical calculi.
2. *Functional or semantic* logical calculi.
3. *Model or pragmatic* logical calculi.

Here we are mostly interested in syntactic logical calculi. In the majority of cases, such calculi are logical system used to prove true formulas (called theorems) and model argumentation and reasoning. Basic building blocks for syntactic logical calculi are formulas from logical languages, algorithms/rules of logical inference and formula transformations, and metarules that determine how to apply algorithms/rules of logical inference to elements from logical languages.

Let $L$ be a logical language or a language of well-formed formulas and $R$ be an algorithmic language, procedural language or a language of rules of inference in $L$. All expressions from $R$ are descriptions of algorithms that work with words from $L$. A standard example of $L$ is the first order predicate language or any language of mathematical logics. But practically $L$ may be any language: natural, mathematical, programming, of chemical formulas, etc. Any language in which it is possible to describe inference rules is an example of $R$. However, $R$ usually contains not only inference rules but also transformation rules.

**Definition 3.1.** A *syntactic* (*deductive*) *logical calculus* in the pair of languages $(L, R)$ is a triad (a named set (Burgin, 2004b)) of the form $\boldsymbol{C} = (A, H, T)$ where $H \subseteq R$, $A, T \subseteq L$ and $T$ is obtained by applying algorithms/procedures/rules from $H$ to elements from $A$.

We consider two main types of deductive logical calculi: exact and relaxed (soft).

In an *exact deductive logical calculus* $\boldsymbol{C} = (A, H, T)$, the system $H$ consists of algorithms.

In an *relaxed deductive logical calculus* $\boldsymbol{C} = (A, H, T)$, the system $H$ consists of procedures some of which are not algorithms.

Examples of exact deductive logical calculi are the syllogistic logic, classical propositional/sentential and predicate logics, intuitionistic logic, inductive logic, deontic logic, weak and strong paraconsistent logics, the logic of imperatives, different temporal logics, class logic, relevant logics, a variety of modal logics, the logic of values, the logic of norms, epistemic logic, erotetic logics, and so on and so forth.

Examples of relaxed deductive logical calculi are dialectic logic, transcendental logic (Husserl, 1929), a logic of diagnosis (Tarasov, *et al*, 1989), and a logic of goal control (Ladenko and Tulchinsky, 1988).

Usually, two cases of logical calculi are studied and constructed: monotonic and non-monotonic.

In the case of a monotonic logical calculus $C$, we have $T = \bigcup_{n=1}^{\infty} T_n$ where $T_1$ is equal to $A$, $T_n$ is equal to $H(T_{n-1})$, and $H(M) = \{ r(N); r \in H \text{ and } N \subseteq M \}$. Another way to represent $T$ is to consider the closure $H^*$ of $H$ with respect to compositions of algorithms. Then $T = H^*(A)$.

In the case of a non-monotonic logical calculus $C$, axioms are changing with time and we have $A = \{ A_n ; n = 1, 2, 3, \ldots \}$ and $T = \{ T_n ; n = 1, 2, 3, \ldots \}$ where $A_n$ is the system of axioms at time (period) $n$ and $T_n$ is the system of axioms at the same time. Sometimes rules of inference $H$ also change with time.

When $L$ is a logical language and $H$ consists of rules of logical deduction, $C$ is a deductive calculus. The same syntactic logical calculus can be considered in different pairs of languages $(L, R)$ as the following simple property shows.

**Lemma 3.1.** If $L_1 \subseteq L_2$ and $R_1 \subseteq R_2$, then any syntactic logical calculus in $(L_1, R_1)$ is a syntactic logical calculus in $(L_2, R_2)$.

Usually we do not explicitly indicate in what pair of languages $(L, R)$ a syntactic logical calculus is considered.

**Remark 3.1.** It is possible to build/define a general syntactic calculus taking any language $L$ and algorithmic language $R$ with algorithms that work with words from $L$. Practically, $L$ may be any language: natural, mathematical, programming, of chemical formulas, etc. For instance, when $L$ contains descriptions and denotations of real/complex numbers and functions, while $H$ consists of rules of differentiation/integration, $C$ is the differentiation/integration calculus. Another example of a general syntactic calculus is any universal algebra (e.g., a group, ring or linear algebra over a field of real numbers).

One more example of a general syntactic calculus is a productive representation of a language, i.e., a language is usually represented as a calculus (cf., Definition 2.4).

**Remark 3.2.** In what follows, it is always assumed that $H$ forms a cumulative system, i.e., algorithms from $H$ only add formulas or change them but never exclude formulas. At the same time, it is possible to consider both inclusive and exclusive algorithms. This brings us to nonmonotonous calculi, which formalize nonmonotonous reasoning and nonmonotonous logics.

**Definition 3.2.** a) $A$ is called the *axiom system* (*base* or *generating system* if $C$ is not a logical calculus) of the calculus $C$.

b) $H$ is called the *system of inference rules* of the calculus $C$.

c) $T$ is called the *body* (the set of theorems or set of the deducible expressions) of the calculus $C$. It is constructed by applying algorithms from $H$ to expressions from $A$.

Components of a syntactic calculus $C = (A, H, T)$ are denoted as follows: $A = A(C)$, $H = H(C)$, and $T = T(C)$.

It is possible that the set of axioms $A$ is empty or countably infinite. In the second case, this set is represented by axiom schemata, which contain variables for elements from the corresponding language $L$ or/and from the alphabet of this language.

*Example* **3.1.** Let us consider the classical propositional or sentential logic/calculus. Many systems of propositional calculus have been devised to achieve consistency, completeness, and independence of axioms. All these systems are logically equivalent in the sense of Definition 3.4a. Thus, it is more correct to call these systems not the same calculus but logically equivalent representations of the classical propositional or sentential logic/calculus.

For instance, Kleene (2002) suggests the following list of axioms (axiom schemas) of the classical propositional calculus:

$$\varphi \to (\chi \to \varphi) \tag{1}$$

$$(\varphi \to (\chi \to \psi)) \to ((\varphi \to \chi) \to (\varphi \to \psi)) \tag{2}$$

$$\varphi \to (\chi \to (\varphi \wedge \chi)) \tag{3}$$

$$\varphi \to \varphi \vee \chi \tag{4}$$

$$\chi \to \varphi \vee \chi \qquad (5)$$

$$\varphi \wedge \chi \to \varphi \qquad (6)$$

$$\varphi \wedge \chi \to \chi \qquad (7)$$

$$(\varphi \to \psi) \to ((\chi \to \psi) \to (\varphi \vee \chi \to \psi)) \qquad (8)$$

$$(\varphi \to \chi) \to ((\varphi \to \neg \chi) \to \neg \varphi) \qquad (9)$$

$$\neg \neg \varphi \to \varphi \qquad (10)$$

Usually the system of inference rules has only one rule called *modus ponens*:

$$\varphi, \varphi \to \psi \vdash \psi$$

or the natural/programming language notation

**If φ and φ → ψ, then ψ.**

Other rules are derived from *modus ponens* and then used in formal proofs to make proofs shorter and more understandable. These rules serve to directly introduce or eliminate connectives, e.g., "**If φ and χ, then φ ∧ χ**" ( or φ, χ ⊢ φ ∧ χ ) or "**If φ, then φ ∨ χ**" ( or φ ⊢ φ ∨ χ ).

A standard transformation rule is substitution. This rule is necessary because axiom schemas demand substitution to become axioms and be applied. Namely, formulas (1) – (10) are axioms when the system of inference rules includes the substitution rule (cf. Example 3.2) and are axiom schemas when the system of inference rules has only Modus Ponens (Smullyan, 1962).

*Example* **3.2.** There are other representations of the classical propositional logic/calculus. Let us consider some of them.

Thus, Church (2002) suggests two logically equivalent representations $\mathbf{P}_1$ and $\mathbf{P}_2$ for the classical propositional calculus. The alphabet of $\mathbf{P}_1$ contains: one symbol of the logical operator ⊃, two punctuation symbols [ and ] , a constant symbol *f*, and a countable set of variables *q, p, s, …* .

The system $\mathbf{P}_1$ has the following list of axiom schemata:

$$[p \supset [q \supset p]] \tag{1}$$
$$[[s \supset [p \supset q]] \supset [[s \supset p] \supset [s \supset q]]] \tag{2}$$
$$[[[p \supset f] \supset f] \supset p] \tag{3}$$

The system of inference rules of $\mathbf{P}_1$ contains two elements:

1. *Modus ponens*: $p$, $p \to q$ imply $q$

2. Substitution Rule: $p$ implies $\mathbf{S}^x_q p$ where $\mathbf{S}^x_q$ denotes the substitution of a variable $x$ by a formula $q$.

The alphabet of $\mathbf{P}_2$ contains: two symbols of logical operators $\supset$ and $\sim$, two punctuation symbols [ and ], and a countable set of variables $q, p, s, \ldots$ .

The system $\mathbf{P}_2$ has the following list of axiom schemata:

$$[p \supset [q \supset p]] \tag{1}$$
$$[[s \supset [p \supset q]] \supset [[s \supset p] \supset [s \supset q]]] \tag{2}$$
$$[[\sim p \supset \sim q] \supset [q \supset p]] \tag{3}$$

The system of inference rules of $\mathbf{P}_2$ is the same as the system of inference rules of $\mathbf{P}_1$:

1. *Modus ponens*: $p$, $p \to q$ imply $q$

2. Substitution Rule: $p$ implies $\mathbf{S}^x_q p$ where $\mathbf{S}^x_q$ denotes the substitution of a variable $x$ by a formula $q$.

*Example* **3.3.** Let us consider the classical first-order predicate logic/calculus. Various systems of first-order predicate calculus have been devised to achieve consistency, flexibility, and independence of axioms. All these systems are logically equivalent in the sense of Definition 3.4a.

Shoenfield (2001) suggests the following list of axiom schemata of the classical first-order predicate calculus:

**Propositional Axiom**: $\varphi \lor \neg \varphi$

**Identity Axiom**: $x = x$

**Substitution Axiom**: $\varphi_x[a] \to \exists x\, \varphi$

**Equality Axioms**: a) If $f$ is a symbol of an $n$-ary function from **F**, then

$$x_1 = y_1 \wedge x_2 = y_2, \ldots, x_n = y_n \rightarrow f(x_1, x_2, \ldots, x_n) = f(y_1, y_2, \ldots, y_n);$$

b) If $f$ is a symbol of an $n$-ary predicate from **P**, then

$$x_1 = y_1 \wedge x_2 = y_2, \ldots, x_n = y_n \rightarrow p(x_1, x_2, \ldots, x_n) = p(y_1, y_2, \ldots, y_n)$$

The system of inference rules of the classical first-order predicate calculus contains two elements:

**Extension rule**: $\varphi$ implies $\varphi \vee \psi$

**Cancellation rule**: $\varphi \vee \varphi$ implies $\varphi$

**Associative rule**: $(\varphi \vee \psi) \vee \chi = \varphi \vee (\psi \vee \chi)$

**Cut rule**: $(\varphi \vee \psi)$ and $(\neg \varphi \vee \chi)$ imply $(\psi \vee \chi)$

**∃-introduction rule**: $\varphi \rightarrow \psi$ implies $\exists x\, \varphi \rightarrow \psi$ if $x$ is not a free variable in $\psi$.

*Example* **3.4.** Syntactic logical calculi provide functional formalization to the notion of a formal theory (Smullyan, 1962). In turn, a formal theory formalizes some source theory from a scientific discipline (e.g., mathematics, physics or economics). In order to specify a formal theory, one first chooses a small collection of predicates, functions and relations, which are regarded as basic for a given field of study (groups, topological spaces or geometry). The chosen predicates delimit the scope of the formal theory. These predicates are the primitives of the theory and together with logical and punctuation symbols (such as the symbols **(** or **,** ) form the alphabet of the theory (calculus) language. The language consists of expressions (functions, relations, and predicates) defined in terms of the primitives. Using them, one writes down certain predicates that are regarded as basic or self-evident within the given field of study. These predicates are the axioms of the theory. It is crucial to make all of underlying assumptions of the source theory explicit as axioms. Often this is not a simple task. One can compare formalization of the Euclidean geometry given by Euclid and Hilbert (Hilbert, 1930). Using logical rules/algorithms of inference (usually, it is only Modus Ponens and substitution rule) theorems of the theory are deduced from the axioms. As a result, a *formal theory* is this structure of theory language, axioms, and theorems.

The process of codifying a scientific discipline by means of primitives and axioms in the predicate calculus is known as formalization. The key issue here is the choice of primitives and axioms. They can be chosen arbitrarily but it is better to exercise a certain aesthetic touch and use the following principles: it must not be too many axioms; they must be basic and self-evident from the discipline's point of view; and they must account for the largest possible number of other concepts and facts.

In all given examples, the system $H$ consists only of simple rules such as *modus ponens* and substitution. Thus it is possible to ask a question why in the definition of a calculus, it is necessary to consider algorithms and not only of simple rules such as *modus ponens* and substitution. The following example explains such a necessity.

*Example* **3.5.** Applications of logic, such as program verification, demand utilization of decidable inference rules. In logic/calculus with validation LV, all inference rules are decidable. In this calculus, $H$ contains a variety of different algorithms in addition to inference rules. In particular, we have a decidable *modus ponens*:

$$\text{If } p \to q \text{ and it is validated that } p \text{ is true, then } q \text{ is true.}$$

To validate that a proposition is true may have different meanings. It can mean:

a) To prove that $p$ is true in the classical sense.
b) To prove that $p$ is true with ordinal induction.
c) To test that $p$ is true.

In such a system, the pair of propositions $p \to q$ and $p$ do not imply $q$ if there is no algorithms in $H$ to check $p$. If validation is temporary, i.e., it true only for some time, then the logic/calculus with validation LV is nonmonotonic. Examples from mathematics of such temporary validation are given in the book of Lakatos (1976).

Other kinds of inference rules (algorithms) are used in fuzzy logics (Zadeh, 1975; Zimmermann, 1997).

**Proposition 3.1.** If $C = (A, H, T)$ and $D = (B, K, Q)$ are syntactic logical calculi in $(L, R)$, $A \subseteq B$, and $H \subseteq K$, then $T \subseteq Q$.

Let us assume that the algorithmic language $R$ contains an identity algorithm $E$, for which $E(w) = w$ for any expression $w$ from $L$.

**Lemma 3.2.** Any subset $Q$ of $L$ is the body of some syntactic logical calculus $C = (A, H, Q)$.

This result shows that the concept of a syntactic logical calculus is very general and to get models better suited to tasks of logic, it is necessary to have some restrictions on those algorithms that are used in logical calculi.

**Lemma 3.3.** For any syntactic logical calculus $C = (A, H, T)$, we have $A \subseteq T$.

**Definition 3.3.** If $A = L$, then a syntactic logical calculus $C = (L, H, T)$ is called a *free syntactic logical calculus* or a *formal deduction system*.

Let us assume that $H$ contains an identity algorithm $E$.

**Corollary 3.1.** The body of a free syntactic logical calculus $C = (L, H, T)$ is equal to $L$.

**Remark 3.3.** There are free syntactic logical calculi in which the language $L$ is infinite, but the body is finite. For instance, we can take the language of the classical propositional calculus as $L$ and such rules of classical deduction that work only with formulas the length of which is less than 1000 as set of inference rules $H$. In this case, the body $T$ of the calculus $C = (L, H, T)$ is finite.

Let $L_1$ and $L_2$ be logical languages, $R_1$ and $R_2$ be algorithmic languages with algorithms that work with words from $L_1$ and $L_2$, correspondingly, and $f: L_1 \to L_2$ is a one-to-one mapping (bijection) of $L_1$ onto $L_2$.

**Definition 3.4.** a) Two syntactic logical calculi $C = (A, H, T)$ in $(L_1, R_1)$ and $D = (B, K, Q)$ in $(L_2, R_2)$ are called *logically* (or *semantically*) *equivalent with respect to $f$* if $Q = f(T)$. This relation is denoted by $C \sim^f_{\lg} D$.

b) Two syntactic logical calculi $C = (A, H, T)$ in $(L, R)$ and $D = (B, K, Q)$ in $(L, R)$ are called *logically* (or *semantically*) *equivalent* if $T = Q$. It is denoted by $C \sim_{\lg} D$.

**Proposition 3.2.** Relation $\sim^f_{\lg}$ is an equivalence relation if and only if $L_1 = L_2$ and $f$ is the identity mapping.

Indeed, when $f$ is not the identity mapping, we have $f \neq g = f^{-1}$ and thus, the relation $C \sim^f_{lg} D$ does not imply the relation $C \sim^g_{lg} D$. Consequently, $\sim^f_{lg}$ is not reflexive and thus, not an equivalence relation.

Sufficiency of conditions in Proposition 3.2 follows from Definition 3.4.

**Corollary 3.2.** Relation $\sim_{lg}$ is an equivalence relation.

Axioms (or their schemata) and rules of inference define a proof theory. Such proof theories are usually considered equivalent when the corresponding calculi are logically equivalent. For instance, various equivalent proof theories of propositional calculus have been constructed (cf. Examples 3.1 and 3.2).

**Definition 3.5.** a) Two syntactic logical calculi $C = (A, H, T)$ in $(L_1, R_1)$ and $D = (B, K, Q)$ in $(L_2, R_2)$ are called *algorithmically* (or *inferentially*) *equivalent* with respect to $f$ if $B = f(A)$ and $Q = f(T)$. It is denoted by $C \sim^f_{alg} D$.

b) Two syntactic logical calculi $C = (A, H, T)$ in $(L, R)$ and $D = (B, K, Q)$ in $(L, R)$ are called *algorithmically* (or *inferentially*) *equivalent* if $T = Q$ and $A = B$. It is denoted by $C \sim_{alg} D$.

Algorithmically equivalent logical calculi are usually considered as different axiomatizations of the same logic.

**Proposition 3.3.** Relation $\sim^f_{alg}$ is an equivalence relation if and only if $L_1 = L_2$ and $f$ is the identity mapping.

**Corollary 3.3.** Relation $\sim_{alg}$ is an equivalence relation.

**Lemma 3.4.** Algorithmic equivalence (with respect to $f$) implies logical equivalence (with respect to $f$).

Syntactic logical calculi can be named by classes of algorithms to which their inference rules belong. For instance, there are finite automaton calculi, recursive calculi, and superrecursive calculi.

**Proposition 3.4.** If the body $T$ of a syntactic logical calculus $C = (A, H, T)$ is finite, then $C$ is logically equivalent to a nondeterministic finite automaton calculus.

*Proof.* Let us consider a syntactic logical calculus $C = (A, H, T)$ with the finite set $T$ of theorems. Taking this finite set $T$ of formulas and choosing some formula $w$ from $T$, it is possible to build a deterministic finite automaton $A_w$ with ε-transition that given the empty word ε as its input, computes $w$ and nothing else. When some other symbol is given to $A_w$, its output is ε.

The finite automaton $A_T$, which computes the set $T$, contains all automata $A_w$ with $w \in T$. It is possible because there are only finite number of such automata $A_w$. The automaton $A_T$ has the start state $q_0$ that is different from the start states of all $A_w$ and works according to the following rules.

With the empty input, the automaton $A_T$ makes a transition to the start state $q_{0w}$ of one of the automata $A_w$. As $A_T$ is a nondeterministic automaton, it has a possibility to make a transition from $q_0$ to any of the states $q_{0w}$. After this transition, the automaton $A_w$ computes the formula $w$, which is produced as the output of $A_T$. In such a way, the automaton $A_T$ computes (deduces) all elements from $T$.

When some other symbol is given to $A_T$, it gives no output as there are no transition from $q_0$ beside ε-transitions. Consequently, if we take $K = \{ A_T \}$, we obtain a nondeterministic finite automaton calculus $C_0 = (A, K, T)$, which is algorithmically and logically equivalent to $C = (A, H, T)$.

Proposition is proved.

**Remark 3.4.** It is possible to prove by the same technique that $C = (A, H, T)$ is algorithmically and logically equivalent to a deterministic finite automaton with ε-transitions syntactic calculus.

**Remark 3.5.** If ε-transitions are not permitted, then the result of Proposition 3.4 is not true for deterministic finite automaton syntactic calculi. For example, when the set of axioms $A$ is empty, the set of theorems a deterministic finite automaton syntactic calculus is also empty.

**Definition 3.6.** Two syntactic logical calculi $C = (A, H, T)$ and $D = (B, K, Q)$ are called *axiomatically* (or *generatively*) *equivalent* if $T = Q$ and $H = K$. It is denoted by $C \sim_{ax} D$.

Axiomatic equivalence of calculi informally means that given the same inference rules/algorithms, different axiom systems produce the same set of theorems.

Definitions imply the following results.

**Lemma 3.5.** The relation $\sim_{ax}$ is an equivalence relation.

**Lemma 3.6.** Axiomatic equivalence implies logical equivalence.

**Lemma 3.7.** Two syntactic logical calculi $C = (A, H, T)$ and $D = (B, K, Q)$ are axiomatically and inferentially equivalent if and only if they coincide.

Let us consider some specific classes of syntactic logical calculi.

**Definition 3.7.** A syntactic calculus $C = (A, H, T)$ is called:

(1) *constructing* if by means of the algorithms (rules) from $H$ new constructions are elaborated;

(2) *transforming* if by means of the algorithms (rules) from $H$ expressions from $L$ are only transformed;

(3) *closed with respect to A* if all elements from $A$ are used for the construction of the set $T$;

(4) *closed with respect to H* if all algorithms from $H$ are used for the construction of the set $T$;

(5) *basically closed* if any algorithm from $H$ may be applied to any set of expressions from $A$;

(6) *transitively closed* if any sequential composition of algorithms from $H$ is admissible;

(7) *completely closed* if it is closed with respect to $A$ and $H$ and is transitively closed;

(8) *admissible* if $T \neq L$, i.e., the set of theorem does not coincide with the whole language;

(9) *consistent with a subset P* of the language *L* if $T \cap P = \emptyset$ where $\emptyset$ is an empty set;

(10) *complete with respect to a set F* of algorithms/rules from *R and a subset Q* of *L* if the set *Q* can be constructed from *T* by means of *F*, i.e. $Q \subseteq F(T)$;

(11) *consistent* if *T* does not include contradictory (false) expressions from *L* (for example, expressions having the form $a \& \neg a$ when *L* contains standard logical connectives);

(12) *complete with respect to a mapping f*: $L \to L$ if for any $a \in L$ we have $a \in T$ or $f(a) \in T$.

For example, when *L* is a logical language with negation $\neg$ and we take $f(a) = \neg a$ then completeness with respect to such *f* is the conventional completeness of a logical calculus. If we take any classical logical calculus $\mathbf{C} = (A, H, T)$ and *P* consists of some expressions having the form $a \& \neg a$, then consistency with *P* is equivalent to conventional definitions of consistency. From the point of view of classical logic, we are compelled to derive any conclusion from inconsistent premises. Consequently, any admissible classical calculus is consistent.

At the same time, when we have an admissible but inconsistent knowledge system, accumulating reports of empirical observations can help in deciding in favor of one alternative over another, allowing one to restore consistency in many cases. However, practical situations show that even without restoring consistency, an inconsistent system can still produce useful information. Examples of admissible but not necessarily consistent syntactic logical calculi are given by different paraconsistent logics (da Costa, 1974; Arruda, 1980). Relevant logics (Anderson and Belnap, 1975) give another class of examples of admissible but possibly inconsistent syntactic logical calculi. In these logics, efficient proof procedures infer only "relevant" conclusions with varying degrees of accessibility, as stated by the criteria of non-classical relevant entailment.

Many of properties from Definition 3.7 are the same for logically equivalent calculi.

**Proposition 3.5.** If $\mathbf{C} = (A, H, T)$ is an admissible (consistent with a subset *P* of the language *L*, complete with respect to a set *F* of algorithms/rules from *R* and a subset *Q* of

$L$, consistent, complete with respect to a mapping $f: L \to L$) syntactic logical calculus and a syntactic logical calculus $\mathbf{B} = (B, K, Q)$ is logically equivalent to $\mathbf{C}$, then $\mathbf{B}$ is an admissible (consistent with a subset $P$ of the language $L$, complete with respect to a set $F$ of algorithms/rules from $R$ and a subset $Q$ of $L$, consistent, complete with respect to a mapping $f: L \to L$, respectively) syntactic logical calculus.

Let us assume that all calculi are completely closed.

**Proposition 3.6.** a) If $T \subseteq L$, $H \subseteq R$ and $H(T) = T$, then $T$ is the body of some transitively closed calculus with the system of inference rules $H$. b) If $T$ is the body of a transitively closed calculus $\mathbf{C} = (A, H, T)$ and $A \subseteq T$, then $H(T) = T$.

Let $H$ includes an algorithm that defines the identity function on $L$. Then Proposition 3.2 implies the following results.

**Corollary 3.4.** If $H(T) \subseteq T$, then $T$ is the body of some transitively closed calculus with the inference rules $H$.

**Corollary 3.5.** The intersection of any set of bodies of transitively closed calculi is the body of some transitively closed calculus.

Any ordinary logical calculus will be completely closed if to the inference rules we add the identity operator on $L$.

**Proposition 3.7.** A calculus $\mathbf{C} = (A, H, T)$ is transitively closed if and only if $H(T) \subseteq T$.

Let $P \subseteq L$ and the whole $L$ may be constructed by applying algorithms from $H$ to an arbitrary element p from $P$.

**Proposition 3.8.** A transitively closed calculus $\mathbf{C} = (A, H, T)$ is admissible if and only if it is consistent with $P$.

Usually relations of some arity are considered in a set as subsets of direct powers of this set (cf. Denotations and Basic Definitions). In logic, we need a more general kind of relations defined between subsets of formulas and individual formulas. Let $X$ be a set. Then (Béziau, 1994), a subset $R$ of the direct product $P(X) \times X$ is called an *abstract logic* or *logical structure*.

The relation *R* is an abstract form of a consequence, or inference, relation in *X*. For instance, the *inference relation* ⊢ in conventional logic is derived from the *consequence relation* (denoted by ⊩ ) of this logic. This consequence relation is called operational. The consequence relation is a subset of the inference relation. Application of one deduction rule, such as *modus ponens*, gives the consequence relation. Consecutive application of several deduction rules gives the inference relation. However, when the set of algorithms *H* is closed under sequential composition, inference and consequence relations coincide. Thus, the inference relation is a kind of the consequence relation.

Besides, the *operational consequence relation* ⊩ there are other types of consequence relations. For instance, there is the *model consequence relation* ⊨ when $\Phi \Vdash \psi$ if $\psi$ is true in all models in which $\Phi$ is true.

In non-monotonic logics, the operational consequence relation $\mid\sim$ called a conditional assertion is used (Kraus, et al, 1990). This relation $\Phi \mid\sim \psi$ is interpreted as *if $\Phi$, normally $\psi$*, or *$\psi$ is a plausible consequence of $\Phi$*.

The consequence relation in an abstract logic determines a formal logical syntax. Logical syntax is not the same as the syntax of a logical language. The syntax or, more exactly, linguistic syntax of a logical language is defined/given by the rules of well-formed formulas construction.

**Definition 3.8.** A consequence relation ⊩ [in a calculus *C* = (*A*, *H*, *T*) ] is defined by the inference system *H* in the following way: for any system of formulas $\Phi$ and formula $\psi$, we have $\Phi \Vdash \psi$ if $\psi \in H(\Phi)$ [$\Phi \Vdash_C \psi$ if $\psi \in H(A \cup \Phi)$ ].

In this case, *H* is called an *operational semantics* for the consequence relation ⊩ .

**Definition 3.9.** A inference relation ⊢ [in a calculus *C* = (*A*, *H*, *T*) ] is defined by the inference system *H* in the following way: for any system of formulas $\Phi$ and formula $\psi$, we have $\Phi \vdash \psi$ if $\psi \in H^*(\Phi)$ [$\Phi \vdash_C \psi$ if $\psi \in H^*(A \cup \Phi)$ ].

An operational semantics for the consequence relation $\Vdash_C$ is called an operational semantics for the logical calculus *C*. The operational semantics for the logical calculus *C* can consists of algorithms of a definite type, e.g., recursive algorithms or inductive algorithms (Burgin, 2005). The type of algorithms determines the type of operational semantics. For example, all conventional logics have a *recursive operational semantics*. *Inductive operational semantics* allows one to use inductive Turing machines for logical inference.

Usually, logic is built from a logical calculus by adding *truth* (functional) and/or *model* (representational) *semantics*. The classical truth semantics is a mapping of all well-formed formulas of a logical language *L* into the set {True, False} or {T, F} or {1, 0}. It is always supposed all axioms, i.e., elements from the set *A*, are true. The truth semantics of a fuzzy logic is a mapping of all well-formed formulas of a logical language *L* into the set [0, 1]. The truth semantics of an intuitionistic logic is a mapping of all well-formed formulas of a logical language *L* into the set {True, False, Unknown}.

In this context, the first Gödel incompleteness theorem has the following form.

**Theorem 3.1.** The truth semantics of a system that contains the formal arithmetic cannot have (be defined by) a recursive operational semantics.

A logical calculus with an inductive operational semantics is more powerful than a logical calculus with a recursive operational semantics and the same system of axioms. For instance, the main result from (Burgin, 2003) implies the following result.

**Theorem 3.2.** The truth semantics of the formal arithmetic can have (be defined by) an inductive operational semantics.

Usually each step of logical inference or deduction involves only a finite number of formulas. Thus, it is more reasonable to consider only finitely based abstract logics.

**Definition 3.10.** A subset *R* of the direct product $P_{fin}(X) \times X$ is called a *finitely based* abstract logic or logical structure.

**Lemma 3.8.** It is possible to represent any finitely based inference relation $R$ in $X$ as a union of $n$-ary relations in $X$, i.e., $R \cong \cup_{n=1}^{\infty} R_n$ where $(\{x_1, x_2, \ldots, x_k\}, z) \in R$ if and only if $(x_1, x_2, \ldots, x_k, z) \in R_{k+1}$ for all $k = 0, 1, 2, \ldots, n, \ldots$.

In logic, inference relations have additional properties. For instance, Kraus, et al, (1990) write, "Reflexivity (i.e., $p \to p$) seems to be satisfied universally by any kind of reasoning based on notion of consequence." However there are logics in which this is not true.

*Example* **3.6.** Let us consider a prediction logic/calculus PV. Its goal is to deduce future from the past and/or present. In this context, knowing that some proposition/predicate $p$ has been true even for a long time, we cannot assert that $p$ will be true tomorrow.

The most notable example of such a situation is given by the famous model of a true empirical proposition that is attributed to Aristotle:

*All swans are white.*

Europeans had believed in this until they came to Australia where they found black swans and disproved this statement.

Many scientific laws have found themselves in a similar situation. For instance, Newton's laws were considered as absolute truth for centuries. However, relativity theory and quantum mechanics demonstrated limitations of these laws. A similar story happened to the famous Church-Turing Thesis, which one of the cornerstones of the contemporary computer science. For a long time, this Thesis was considered as absolute truth. Only recently with the advent of superrecursive algorithms, it was refuted (Burgin, 2005).

**Definition 3.11.** A formal consequence/inference relation $R$ in $X$ is called $m$-bounded if $R \cong \cup_{n=1}^{m} R_n$.

**Definition 3.12.** A formal consequence/inference relation $R$ in $X$ is called functionally $m$-bounded if the range $\text{Rg}(R) = \cup_{n=1}^{m} \text{Rg}(Q_n)$ where each $Q_n$ is an $n$-ary relation in $X$ for $n = 1, 2, \ldots, m$.

Informally, functional *m*-boundedness means that to get the range of *R*, we need only *m*-ary relations and relations of smaller arity. The range is important because if *R* is an inference relation in some actual logic, its range consists of all theorems of this logic.

**Definition 3.13.** A formal consequence/inference relation *R* in *X* is called *m*-strict if $R \cong R_m$.

**Definition 3.14.** A formal consequence/inference relation *R* in *X* is called functionally *m*-strict if $\mathrm{Rg}(R) = \mathrm{Rg}(Q_m)$ where $Q_m$ is an *m*-ary relation.

Let *H* be a system of (inference) algorithms/rules in a logical language *L*.

**Lemma 3.9.** If the system *H* is closed with respect to the sequential composition of algorithms, then $\vdash = \Vdash$.

**Theorem 3.3.** For any syntactic logical calculus $C = (A, H, T)$, we have $T = \{ \varphi \in L; A \Vdash_C \varphi \}$.

It is useful to specify the concept of a logical structure, making it closer to concrete logics and logical calculi.

Let *L* be an arbitrary language and *A* and *T* be subsets of *L*.

**Definition 3.16.** An abstract universal logical calculus *U* in a set *X* with a base *A* and inference relation *F* is a triad $U = (A, F, T)$ such that $F \subseteq P(X) \times X$ and $T = F(A)$.

Let *L* be a logical language and *A* and *T* be subsets of *L*.

**Definition 3.17.** A formal universal logical calculus is an abstract universal logical calculus in a logical language *L*.

**Theorem 3.4.** The consequence relation $\Vdash_C$ in *C* is functionally *m*-bounded in *C* if and only if for any pair $(\{x_1, x_2, \ldots, x_k\}, x)$ from $\Vdash_C$ with $k > m$ there is a set $\{x_1, x_2, \ldots, x_n\}$ of element *z* in *L* such that $(\{x_1, x_2, \ldots, x_n\}, x)$ also belongs to $\Vdash_C$ and $n < m + 1$.

**Corollary 3.3.** The consequence relation $\Vdash_C$ in $C$ is functionally 2-bounded in $C$ if and only if for any pair ($\{x_1, x_2, \ldots, x_k\}, x$) from $\Vdash_C$ there is an element $z$ in $L$ such that $(z, x)$ also belongs to $\Vdash_C$.

**Theorem 3.5.** The inference relation $\vdash_C$ in $L$ with respect to $C$ is functionally $m$-bounded in $C$ if and only if for any pair ($\{x_1, x_2, \ldots, x_k\}, x$) from $\vdash_C$ with $k > m$ there is a set $\{x_1, x_2, \ldots, x_n\}$ of element $z$ in $L$ such that ($\{x_1, x_2, \ldots, x_n\}, x$) also belongs to $\vdash_C$ and $n < m + 1$.

**Corollary 3.3.** The inference relation $\vdash_C$ in $L$ with respect to $C$ is functionally 2-bounded in $C$ if and only if for any pair ($\{x_1, x_2, \ldots, x_k\}, x$) from $\vdash_C$ there is an element $z$ in $L$ such that $(z, x)$ also belongs to $\vdash_C$.

## 9. Conclusion

Basic properties of such fundamental logical systems as syntactic logical calculi have been explored. Relations between axiom systems (such as axiomatic equivalence), systems of inference rules/algorithms (such as algorithmic equivalence), and different calculi ((such as logical equivalence) have been described. The inherent structure of logical calculus has been explicated and formalized as the first in the development of metalogic. This approach allows one to study how logical languages are built and what is a general schema of theorem derivation from axioms in a logical calculus.

There are several directions for future research. It would be interesting to study the inner structure of the set of theorems. It is possible to stratify this set by complexity of theorem inference, time of inference or by other relations between theorems. Inference complexity can be estimated by the length of inference, by the maximal length of formulas utilized in the inference procedure or by the complexity of algorithms used for inference.

Another important direction is to study other types of calculi, that is, model and functional calculi, from an algorithmic perspective.